\documentclass{amsart}
\usepackage[headings]{fullpage}
\usepackage{epsf, amssymb, amsfonts, amsmath}

\newtheorem{thm}{Theorem}
\newtheorem{conj}[thm]{Conjecture}
\newtheorem*{thm*}{Theorem}
\newtheorem{lem}[thm]{Lemma}
\newtheorem{prop}[thm]{Proposition}
\newtheorem{cor}[thm]{Corollary}
\theoremstyle{definition}

\theoremstyle{remark}

\newcommand{\Z}{{\mathbb{Z}}}

\newcommand{\cen}{{\rm Z}}

\newcommand{\Nalt}{\mathop{\rm N_{alt}}\nolimits}
\newcommand{\Sym}{\mathop{\rm Sym}\nolimits}
\newcommand{\LNalt}{\mathop{\rm LN_{alt}}\nolimits}
\newcommand{\Tder}{\mathop{\rm Tder}\nolimits}
\newcommand{\N}{{\mathbb{N}}}
\newcommand{\ant}{{S}}

\newcommand{\ddx}[1]{{\frac{\partial}{\partial x_{#1}}}}

\def\Nuc{\mathop{\rm N}\nolimits}
\def\Endo{\mathop{\rm End}\nolimits}
\def\ad{\mathop{\rm ad}\nolimits}
\def\alg{\mathop{\rm alg}\nolimits}
\def\spann{\mathop{\rm span}}

\def\Gr{\mathop{\rm Gr}}

\def\L{{\mathcal L}}

\title[Universal enveloping algebras of Lie triple systems]{Ideals in non-associative
universal enveloping algebras of Lie triple systems}
\author{J.Mostovoy \and J.M. P\'erez-Izquierdo}
\address{Instituto de Matem\'aticas, Unidad Cuernavaca,
 Universidad Nacional Aut\'onoma de M\'exico,
A.P. 273-3,  C.P. 62251, Cuernavaca, Morelos, MEXICO}
\email{jacob@matcuer.unam.mx}
\address{Departamento Matem\'aticas y Computaci\'on, Universidad de La
Rioja, Logro\~no, 26004, SPAIN} \email{jm.perez@dmc.unirioja.es}

\subjclass[2000]{20N05, 17D99}

\thanks{Both authors were partially supported by the SEP-CONACyT
grant no.\ 44100. J.M.P.-I.\ acknowledges support from
BFM2001-3239-C03-02 (MCYT) and  ANGI2001/26 (Plan Riojano de I+D).
J.M.\ expresses his gratitude to Max-Planck-Institut f\"ur
Mathematik, Bonn, where part of the work presented here was carried
out.}

\begin{document}

\begin{abstract}
The notion of a non-associative universal enveloping algebra for a
Lie triple system arises when Lie triple systems are considered as
Bol algebras (more generally, Sabinin algebras). In this paper a new
construction for these universal enveloping algebras is given, and
their properties are studied.

It is shown that universal enveloping algebras of Lie triple systems
have surprisingly few ideals. It is conjectured, and the conjecture
is verified on several examples, that the only proper ideal of the
universal enveloping algebra of a simple Lie triple system is the
augmentation ideal.
\end{abstract}

\maketitle


\section{Introduction.}

Given a smooth manifold $M$ and a point $e\in M$, a {\em local
multiplication} on $M$ at $e$ is a smooth map $U\times U\to M$ where
$U$ is some neighbourhood of $e$ and the point $e$ is a two-sided
unit, that is, $xe=ex=x$ for all $x\in U$. If $x$ is sufficiently
close to $e$, both left and right multiplications by $x$  are
one-to-one. Therefore, there always exists a neighbourhood $V\subset
U$ where the operations of left and right division are defined by
the identities $a\backslash (ab)=b$ and $(ab)/b=a$ respectively. Two
local multiplications at the same point $e$ of a manifold $M$ are
considered to be equivalent if they coincide when restricted to some
neighbourhood of $(e,e)$ in $M\times M$. Equivalence classes of
local multiplications are called {\em infinitesimal loops}.
(Sometimes infinitesimal loops are also called {\em local loops}.)

The importance of infinitesimal loops lies in the fact that they are
closely related to affine connections on manifolds. Namely, any
affine connection on $M$ defined in some neighbourhood of $e$
determines a local multiplication at $e$. Conversely, each (not
necessarily associative) local multiplication at $e$ defines an
affine connection on some neighbourhood of $e$; this gives a
one-to-one correspondence between germs of affine connections and
infinitesimal loops. The details can be found, for example, 
in \cite{S}.

Local non-associative multiplications on manifolds can rarely be
extended to global multiplications and, thus, cannot be studied
directly by algebraic means. Nevertheless, any local multiplication
gives rise to an algebraic structure on the tangent space at the
unit element, consisting of an infinite number of multilinear
operations.  Such algebraic structures are known as {\em Sabinin
algebras}; for associative multiplications they specialise to Lie
algebras. Given a Sabinin algebra that satisfies certain convergence
conditions, one can uniquely reconstruct the corresponding analytic
infinitesimal loop. Therefore, Sabinin algebras may be considered as
the principal algebraic tool in studying local multiplications and
local affine connections.

The general theory of Sabinin algebras has so far only been
developed over fields of characteristic 0. From now on we shall
assume that this is the case: {unless stated otherwise, all vector
spaces, algebras etc will be assumed to be defined over a field $F$
of characteristic zero.}

Many general properties of Sabinin algebras are similar to those of
Lie algebras. In particular, any Sabinin algebra $V$ can be realised
as the space of primitive elements of some "non-associative Hopf
algebra" $U(V)$, called the {\em universal enveloping algebra} of
$V$. The operations in $V$ are naturally recovered from the product
in $U(V)$. Just as in the Lie algebra case, the universal enveloping
algebras of Sabinin algebras have Poincar\'e--Birkhoff--Witt bases.
If a Sabinin algebra $V$ happens to be a Lie algebra, $U(V)$ is
precisely the usual universal enveloping algebra of the Lie algebra
$V$.

The definition of a Sabinin algebra involves an infinite number of
multilinear operations that satisfy rather complicated identities;
we refer to \cite{P2}, \cite{SM-Orange} or \cite{SU} for the precise
form of these. However, additional conditions imposed on a local
multiplication may greatly simplify the structure of the
corresponding Sabinin algebra. For example, the associativity
condition implies that only one of all the multilinear operations is
non-zero; the identities of a Sabinin algebra specialise to the
identities defining a Lie algebra, that is, antisymmetry and the
Jacobi identity. If a local multiplication satisfies the Moufang law
\begin{equation}
\label{Moufang.eq} a(b(ac)) = ((ab)a)c \quad \hbox{and} \quad
((ca)b)a = c(a(ba)),
\end{equation}
the corresponding Sabinin algebra is a {\em Malcev algebra}. A
vector space with a bilinear skew-symmetric operation (bracket) is
called a Malcev algebra if the bracket satisfies
\begin{equation*}
[J(a,b,c),a] = J(a,b,[a,c])
\end{equation*}
where $J(a,b,c) = [[a,b],c] + [[b,c],a] + [[c,a],b]$ denotes the
jacobian of $a,b$ and $c$.

Imposing the left Bol identity
\begin{equation*}
a(b(ac)) = (a(ba))c,
\end{equation*}
on the local multiplication, we obtain the structure of a {\em left
Bol algebra} on the tangent space to the unit. A left Bol algebra is
a vector space with one bilinear and one trilinear operation,
denoted by $[\,,\,]$ and $[\,,\,,\,]$ respectively. The ternary
bracket must satisfy the following relations:
\begin{eqnarray*}
&[a,a,b] = 0 & \\
&[a,b,c] + [b,c,a] + [c,a,b] = 0 &\\
&[x,y,[a,b,c]] =[[x,y,a],b,c] + [a,[x,y,b],c] + [a,b,[x,y,c]].
\end{eqnarray*}
The binary bracket is required to be skew-symmetric and should
satisfy
\begin{equation*}
[a,b,[x,y]] = [[a,b,x],y] + [x,[a,b,y]] +[x,y,[a,b]] +
[[a,b],[x,y]].
\end{equation*}

Bol algebras generalise Malcev algebras. Indeed, in any Malcev
algebra a ternary bracket can be defined by
\begin{equation*}
[a,b,c] = [[a,b],c] - \frac{1}{3} J(a,b,c).
\end{equation*}
With this additional operation a Malcev algebra becomes a Bol
algebra.

Another important subclass of Bol algebras are {\em Lie triple
systems}; these are the Bol algebras whose binary bracket is
identically equal to zero. Lie triple systems arise as tangent
spaces to smooth local {\em Bruck loops} (also known as {\em
K-loops}). These loops, in addition to the left Bol identity,
satisfy the identity
\begin{equation*}
(ab)^{-1} = a^{-1} b^{-1}.
\end{equation*}
where $x^{-1}$ is shorthand for $e/x$; see \cite{Ki}. Lie triple
systems play a prominent role in the theory of symmetric spaces
since a symmetric space can be given the structure of a local Bruck
loop at any point.

\medskip

Identities satisfied in an infinitesimal loop can be translated into
identities satisfied in the universal enveloping algebra of the
corresponding Sabinin algebra. In particular, the universal
enveloping algebra $U(M)$ of a Malcev algebra $M$ is a
non-associative bialgebra that satisfies the linearisations
\begin{equation*}
\sum a_{(1)}(y(a_{(2)}z)) = \sum ((a_{(1)} y)a_{(2)}) z
\end{equation*}
and
\begin{equation*}
\sum ((ya_{(1)})z)a_{(2)}= \sum y(a_{(1)} (z a_{(2)}))
\end{equation*}
of (\ref{Moufang.eq}). Here we use Sweedler's notation \cite{Sw} for
the comultiplication: $\Delta(a) = \sum a_{(1)}\otimes a_{(2)}$.
Since $M$ coincides with the subspace of all primitive elements of
$U(M)$ we have $\Delta(a) = a\otimes 1 + 1\otimes a$ for any $a\in
M$, and, hence, $a(yz) + y(az) = (ay)z + (ya)z$ and $(ya)z + (yz)a =
y(az) + y(za)$, or, equivalently
\begin{equation}
\label{alternative.eq} (a,y,z) = -(y,a,z) = (y,z,a).
\end{equation}
Therefore, $M$ lies in the generalised alternative nucleus
$\Nalt(U(M))$ of $U(M)$. (The subset $\Nalt(A)$ of an algebra $A$
consists of all $a\in A$ that satisfy (\ref{alternative.eq}) for any
$y,z\in A$). The product on $M$ is recovered as $[a,b] = ab - ba$ in
$U(M)$.

The universal enveloping algebra $U(V)$ of a Bol algebra $V$
satisfies the identity
\begin{equation}
\label{BolHopf.eq} \sum a_{(1)}(y(a_{(2)} z)) = \sum
(a_{(1)}(ya_{(2)}))z.
\end{equation}
Since $V$ coincides with the primitive elements of $U(V)$, for any
$a\in V$ and $y,z \in U(V)$ we have that
\begin{equation}
\label{leftalternative.eq} (a,y,z) = -(y,a,z).
\end{equation}
This is equivalent to saying that $V$ is contained in the left
generalised alternative nucleus $\LNalt(U(V))$ of the algebra
$U(V)$. The binary and the ternary products on $V$ are recovered by
\begin{equation*}
[a,b] = ab - ba \quad \hbox{and}\quad [a,b,c] = a(bc) - b(ac) -
c(ab) + c(ba)
\end{equation*}
in $U(V)$. It is known \cite{P} that for any algebra $A$
\begin{equation*}
\LNalt(A)=\{ a \in A\,\vert\, (a,x,y) = -(x,a,y) \enskip \forall x,y
\in A\}
\end{equation*}
is a Lie triple system with $[a,b,c] = a(bc) - b(ac) - c(ab) +
c(ba)$.

\medskip

Universal enveloping algebras for Malcev, Bol and general Sabinin
algebras have been introduced only recently; their properties are
still waiting to be explored. It might be tempting to assume that
the theory of universal enveloping algebras for Lie algebras can be
extended rather painlessly to the case of general Sabinin algebras,
especially since many aspects of the theory are known to generalise
well. However, it turns out that some very basic properties, such as
the abundance of ideals in the universal enveloping algebras of Lie
algebras, fail to hold in the general non-associative case. In
particular, we shall see that while the properties of Malcev and Bol
algebras, discussed above, may look similar, this similarity does
not extend too far.

The motivation for this paper is the following version of Ado's
Theorem for Malcev algebras that appeared in \cite{PS}:

\begin{thm}
For any finite--dimensional Malcev algebra $M$ over a field of
characteristic $\neq 2,3$ there exists a unital finite--dimensional
algebra $A$ and a monomorphism of Malcev algebras $\iota \colon M
\rightarrow \Nalt(A)$.
\end{thm}

One is prompted to ask whether a similar statement holds for other
classes of Bol algebras, in particular, for Lie triple systems.
Given a finite dimensional Lie triple system $V$, one could ask
whether it is contained as a subsystem of $\LNalt(A)$, with $ab =
ba$ for all $a,b \in V$, for some finite dimensional unital algebra
$A$. It is easy to see that this happens if and only if there exists
an ideal of finite codimension in $U(V)$ which intersects $V$
trivially. Our answer shows that  for Lie triple systems the
situation is very different from the case of Lie or Malcev algebras:

\begin{thm}\label{main.thm}
Let $A$ be a finite dimensional unital algebra over a field $F$ of
characteristic 0 and $V$ --- a Lie triple system contained as a
subsystem in $\LNalt(A)$ such that $ab= ba$ for all $a,b \in V$.
Assume that $A$ is generated by $V$ as a unital algebra. Then $V$ is
nilpotent and $A$ decomposes (as a vector space) into a direct sum
of a nilpotent ideal and a central subalgebra without nonzero
nilpotent elements.
\end{thm}

Note that we do not claim that embeddings mentioned in
Theorem~\ref{main.thm} do exist for all nilpotent Lie triple
systems.

\medskip

Examples suggest that the ideals in the universal enveloping
algebras of Lie triple systems are even scarcer than it is implied
by Theorem~\ref{main.thm}.
\begin{conj}\label{main.conj}
The only proper ideal of the universal enveloping algebra of a
simple Lie triple system is its augmentation ideal.
\end{conj}
We shall verify the above conjecture in several cases by direct
calculations in Poincar\'e--Birkhoff--Witt bases.

\medskip

For each Lie triple system $V$ there exists a canonically defined
Lie algebra $\L_S(V)$, called the {\em Lie envelope of $V$} of which
$V$ is a subsystem. The Poincar\'e--Birkhoff--Witt Theorem allows to
identify the algebra $U(V)$ with a subspace of $U(\L_S(V))$.
Motivated by analogy with Bruck loops, we shall show how the
multiplication on $U(\L_S(V))$ can be modified to become compatible
with the non-associative multiplication on $U(V)$.

\medskip

The paper is organised as follows. The next section is auxiliary; it
is a loose collection of various properties of Bol algebras and Lie
triple systems. Section~\ref{finite.sec} contains the proof of
Theorem~\ref{main.thm}. The construction of the universal enveloping
algebra of a Lie triple system via its Lie envelope is given in
Section~\ref{construction.sec}. Finally, in Section~\ref{simple.sec}
we present some evidence for Conjecture~\ref{main.conj}.

We have made no attempt to make this paper self-contained. We refer
to \cite{PS} for the properties of the universal enveloping algebras
of Malcev algebras, to \cite{P} --- for Bol algebras and to \cite{P}
--- for general Sabinin algebras. The paper of Lister \cite{Lis} is
the general reference for Lie triple systems; the questions of
nilpotency are treated in \cite{Hop} .

About the notation: we shall often write "L.t.s." for "Lie triple
system". As usual, the true meaning of "non-associative" is "not
necessarily associative"; however "non-nilpotent" stands for "not
nilpotent". The notations $L_x$ and $R_x$ are used to denote the
multiplication by $x$ on the left and on the right respectively; the
sum $L_a+R_a$ is denoted by $T_a$. The product $a(a(\cdots (aa))$
will be written simply as $a^n$. The left, middle and right
associative nuclei of an algebra $A$ are denoted by $\Nuc_l(A)$,
$\Nuc_m(A)$ and $\Nuc_r(A)$ respectively, while $\cen(A)$ is the
notation for the center of $A$. (Recall that the left associative
nucleus of $A$ is the set of all $a\in A$ such that $(a,y,z)=0$ for
arbitrary $y,z\in A$; the right and the middle associative nuclei
are defined similarly.) By $\alg\langle X \rangle$ (or
$\alg_1\langle X \rangle$) we denote the subalgebra (unital
subalgebra, respectively) generated by the subset $X\subset A$.


\section{Some properties of the enveloping algebras for
Bol algebras and Lie triple systems. }\label{stuff.sec}

\begin{lem}
\label{derivation.lem} Let $(V,[\,,\,,\,],[\,,\,])$ be a Bol
algebra. For $a,b\in V$ such that $[a,b] = 0$, the map $[L_a,L_b]$
is a derivation of $U(V)$.
\end{lem}

Recall that a ternary derivation of an algebra $A$ is a triple
$(d_1,d_2,d_3)$ of linear maps such that
\begin{equation*}
d_1(xy) = d_2(x) y + x d_3(y)
\end{equation*}
for all $x,y \in A$. The set $\Tder(A)$ of all ternary derivations
of $A$ is a Lie algebra with the obvious bracket. It is clear that
if $d_1(1) = d_2(1) = d_3(1) = 0$ then $d_1 = d_2 = d_3$ is a
derivation of $A$.

\begin{proof}[Proof of Lemma~\ref{derivation.lem}]
Notice that the identity (\ref{leftalternative.eq}) can be written
as $(L_a,T_a,-L_a) \in \Tder(U(V))$ and, as a consequence,
\[([L_a,L_b],[T_a,T_b],[L_a,L_b]) \in \Tder(U(V)).\]
Evaluating both commutators at $1$, we observe that $[L_a,L_b](1) =
[a,b] = 0 = [T_a,T_b](1)$, so $[L_a,L_b] = [T_a,T_b]$ is a
derivation of $U(V)$.
\end{proof}

\begin{lem}
Let $(V,[\,,\,,\,],[\,,\,])$ be a Bol algebra. For $a,b\in V$ such
that $[a,b] = 0$ and any $x\in U(V)$
\[ [L_a,L_b](x)=-2(a,b,x).\]
\end{lem}
\begin{proof}
The identity (\ref{leftalternative.eq}) with $y = b$ gives  $L_a L_b
+ L_b L_a = L_{ab+ba} = 2L_{ab}$. Therefore, $[L_a,L_b](x) =
L_aL_b(x) -(2L_{ab} - L_aL_b)(x) = -2 (a,b,x)$.
\end{proof}

\begin{lem}\label{leftmult.lem}
Let $(V,[\,,\,,\,],[\,,\,])$ be a Bol algebra. For any $a \in V$
\[L_{a^n} L_{a^m} = L_{a^{n+m}}.\]
\end{lem}
\begin{proof}
See Proposition~38 in \cite{P2}.
\end{proof}

For any Sabinin algebra $V$, the universal enveloping algebra is an
H-bialgebra. That is, $U(V)$ is a non-associative unital bialgebra
equipped with two bilinear maps, $\backslash \colon U(V)\times U(V)
\rightarrow U(V)$ and $/\colon U(V) \times U(V) \rightarrow U(V)$
such that
\begin{eqnarray*}
\sum x_{(1)}\backslash (x_{(2)} y ) &=\qquad \epsilon (x) y\qquad =&
\sum x_{(1)} (x_{(2)} \backslash y ) \hbox{ and }
\\
\sum (y x_{(1)})/ x_{(2)} &=\qquad\epsilon (x) y\qquad =& \sum
(y/x_{(1)}) x_{(2)}.
\end{eqnarray*}
The behaviour of these maps with respect to the comultiplication
$\Delta$ and the counit $\epsilon$ is expressed by
\begin{equation*}
\Delta(x\backslash y) = \sum x_{(1)}\backslash y_{(1)}\otimes
x_{(2)}\backslash y_{(2)}, \quad \Delta(y/x) = \sum
y_{(1)}/x_{(1)}\otimes y_{(2)}/ x_{(2)}
\end{equation*}
and
\begin{equation*}
\epsilon(x\backslash y) = \epsilon(x) \epsilon(y), \quad \epsilon
(y/x) = \epsilon(x) \epsilon(y).
\end{equation*}

Fix an ordered basis $\{a_i\}_{i\in \Lambda}$ of $V$, with $\Lambda$
being the index set.  The algebra $U(V)$ then has the
Poincar\'e--Birkhoff--Witt basis
\[\{ a_{i_1}(a_{i_2}(\cdots
(a_{i_{n-1}} a_{i_n})\cdots))\,\vert\ i_1\leq \cdots \leq i_n \hbox{
and } n \in \N\}.\] The algebra $U(V)$ is filtered by $U(V) =
\cup_{n\in \N} U(V)_n$ with \[U(V)_n=\spann \langle
a_{1}(a_{2}(\cdots (a_{m-1} a_{m})) \,\vert\, a_1,\dots,a_m\in V,
\enskip m\leq n\rangle.\] The {\em degree} of an element of $U(V)$
with respect to this filtration is defined in the obvious way. The
corresponding graded algebra $\Gr{U(V)}$ is isomorphic to $\Sym(V)$,
the symmetric algebra on $V$.

\medskip

Let $(V,[\,,\,,\,])$ be a Lie triple system, and $U(V)$ --- its
universal enveloping algebra. The automorphism $a \mapsto -a$ of $V$
extends to an automorphism $\ant: U(V)\to U(V)$.

\begin{lem}\label{LTScommutators.lem}
Let $(V,[\,,\,,\,])$ be a L.t.s. and $U(V)$ --- its universal
enveloping algebra. Then for any $a \in V$ we have that $[a, U(V)_n]
\subseteq U(V)_{n-1}$.
\end{lem}
\begin{proof}
Let $x = a_{1}(a_{2}(\cdots (a_{{n-1}} a_{n})\dots)) \in U(V)_n$
with $a_1,\dots,a_n \in V$. Since $\Gr{U(V)}$ is isomorphic to
$\Sym(V)$, $[a,x]$ belongs to $U(V)_n$. On the other hand,
$\ant([a,x]) = [-a,(-1)^nx] = (-1)^{n-1} [a,x]$. Therefore, $[a,x]
\in U(V)_{n-1}$.
\end{proof}

The automorphism $\ant$ notably simplifies the left division
$\backslash$ on $U(V)$.

\begin{prop}
\label{Kloop.prop} Let $(V,[\,,\,,\,])$ be a L.t.s. For all $x,y\in
U(V)$
\begin{equation*}
x\backslash y = \ant(x) y \quad \hbox{and}\quad \ant(x) =
x\backslash 1=1/x.
\end{equation*}
\end{prop}
\begin{proof}
Let us prove that $\sum \ant(x_{(1)}) x_{(2)} = \epsilon (x) 1$. To
this end we observe that this is a linear relation, so we only have
to verify it on a set of elements spanning the vector space $U(V)$,
for instance, $\{1\} \cup \{ a^n\,\vert\, a \in V\}$ with $a^n =
a(\cdots (aa))$. We have $\sum \ant({a^n}_{(1)}) {a^n}_{(2)} =
\sum_{k=0}^{n} {n \choose k} \ant(a^k)a^{n-k} = \sum_{k=0}^{n} {n
\choose k} (-1)^k a^n = 0 = \epsilon(a^n)$, as desired.

From (\ref{BolHopf.eq}) and $\sum \ant(x_{(1)}) x_{(2)} = \epsilon
(x) 1$ we obtain
\[\sum x_{(1)}(\ant(x_{(2)})(x_{(3)}y)) = \sum
(x_{(1)}(\ant(x_{(2)})x_{(3)})) y = \sum (x_{(1)} \epsilon(x_{(2)}))
y = xy.\] By the definition of $\backslash$ we have
\begin{equation*}
\sum \ant(x_{(1)})(x_{(2)}y)= \sum x_{(1)}\backslash
(x_{(2)}(\ant(x_{(3)})(x_{(4)}y))) =\sum x_{(1)}\backslash (x_{(2)}
y )= \epsilon(x)y
\end{equation*}
so
\begin{equation*}
\ant(x) y = \sum \ant(x_{(1)})(x_{(2)}(x_{(3)}\backslash y) ) = \sum
\epsilon(x_{(1)})\cdot x_{(2)}\backslash y = x\backslash y.
\end{equation*}
With $y =1$ we get $\ant(x) = x\backslash 1$, and from  $\sum
S(x_{(1)}) x_{(2)} = \epsilon (x) 1$ we also get  $S(x) = \sum
(S(x_{(1)}) x_{(2)})/x_{(3)} = \sum \epsilon (x_{(1)}) 1/x_{(2)} =
1/x$.
\end{proof}

Proposition \ref{Kloop.prop} ensures that $U(V)$ satisfies the
linearisation of the equations defining a Bruck loop. Therefore, the
linearisation of any identity satisfied by Bruck loops will hold in
$U(V)$. Consider, for instance, the so-called {\em precession map}
$\delta_{a,b} \colon c \mapsto (ab) \backslash (a(bc))$. For a Bruck
loop this map is known to be an automorphism \cite{Ki}. Linearising
this result we obtain

\begin{cor}
Let $(V,[\,,\,,\,])$ be a L.t.s. The map $\delta_{x,y} \colon U(V)
\rightarrow U(V)$ given by
\begin{equation*}
\delta_{x,y}(z) = \sum (x_{(1)}y_{(1)})\backslash
(x_{(2)}(y_{(2)}z))
\end{equation*}
satisfies
\begin{equation*}
\delta_{x,y}(wz) = \sum \delta_{x_{(1)},y_{(1)}}(w)
\delta_{x_{(2)},y_{(2)}}(z).
\end{equation*}
\end{cor}
The maps $\delta_{x,y}$ reflect the lack of associativity in $U(V)$.
They satisfy
\begin{equation}
\label{precession.eq} \sum (x_{(1)}y_{(1)})
\delta_{x_{(2)},{y_{(2)}}}(z) = x(yz).
\end{equation}
Clearly, $\Delta(\delta_{x,y}(z)) = \sum
\delta_{x_{(1)},{y_{(1)}}}(z_{(1)})\otimes
\delta_{x_{(2)},{y_{(2)}}}(z_{(2)})$. Thus,
\begin{equation}
\label{zerodegree.eq} \delta_{x,y}(V) \subseteq V
\end{equation}
 and in general
\begin{equation*}
\delta_{x,y}(U(V)_n)\subseteq U(V)_n.
\end{equation*}
The maps $\delta_{x,a}$ and $\delta_{a,x}$ are derivations of $U(V)$
for any $a \in V$. In fact, $\delta_{a,b}(x) = -(a,b,x)$ and
$\delta_{a,b}(c) = \frac{1}{2} [a,b,c]$ for any $a,b,c\in V$ and $x
\in U(V)$.

The following statement is a direct analogue of the corresponding
result for Bruck loops \cite{Ki}.

\begin{prop}
Let $(V,[\,,\,,\,])$ be a L.t.s. Then the left and the middle
associative nuclei of $U(V)$ coincide: \[\Nuc_l(U(V)) =
\Nuc_m(U(V)).\]
\end{prop}
\begin{proof}
The identity (\ref{BolHopf.eq}) implies
\begin{equation}\label{middle.eq}
\sum x_{(1)}\bigl((S(x_{(2)})y) (x_{(3)}z)\bigr) = \sum
\bigl(x_{(1)}((S(x_{(2)})y) x_{(3)})\bigr) z.
\end{equation}
If $y$ is in $\Nuc_m(U(V))$,  the left-hand side of
(\ref{middle.eq}) is equal to
\[\sum x_{(1)}\bigl(S(x_{(2)})(y (x_{(3)}z))\bigr) = y(xz).\]
On the other hand, the right-hand side of (\ref{middle.eq}) can be
re-written as
\[\sum  \bigl(x_{(1)}(S(x_{(2)})(y x_{(3)}))\bigr) z= (yx)z\]
and, hence, $y(xz) = (yx)z$ for all $x, z \in U(V)$. Therefore,
$\Nuc_m(U(V)) \subseteq N_l(U(V))$.

Similarly, notice that (\ref{BolHopf.eq}) also implies
\begin{equation*}
\sum x_{(1)}\bigl((yS(x_{(2)})) (x_{(3)}z)\bigr) = \sum
\bigl(x_{(1)}((yS(x_{(2)})) x_{(3)})\bigr) z.
\end{equation*}
For $y \in \Nuc_l(U(V))$ one concludes that $x(yz) = (xy)z$ for all
$x,z \in U(V)$ and, hence, that $\Nuc_l(U(V)) \subseteq
\Nuc_m(U(V))$.
\end{proof}

\begin{lem}
\label{center.lem} Let $(V,[\,,\,,\,])$ be a L.t.s. and $A$ --- a
quotient of $U(V)$. If $a\in V$ satisfies $[L_a,L_b]=0$ for all
$b\in V$, then $a\in \cen(A)$.
\end{lem}
\begin{proof}
For any $x\in A$ we have $L_x\in\alg_1\langle L_b\,\vert\, b\in
V\rangle$. This can be established by induction on the degree of $x$
with respect to the PBW filtration that $A$ inherits from $U(V)$,
using the fact that $L_{by + yb}=L_b L_y + L_y L_b$ for any $y\in A$
and $b\in V$.

Since $[L_a,L_b]=0$ for all $b\in B$ we have that $[L_a,L_x] = 0$
for any $x\in A$, so $a(xy) =x(ay)$ for any $x,y\in A$. Setting $y
=1$ we get that $ax = xa $ for any $x\in A$. Therefore,
$(xy)a=a(xy)=x(ay)=x(ya)$ and $a\in\Nuc_r(A)$. This can also be
expressed by saying that the triple $(R_a,0,R_a)$ belongs to
$\Tder(A)$.

The identity (\ref{leftalternative.eq}) implies that
$(L_a,T_a,-L_a)$ is also in $\Tder(A)$. Since $R_a =L_a$, it follows
that $(2L_a,2L_a,0)\in \Tder(A)$ and thus $a\in \Nuc_l(A)$.
Similarly, $(0,2R_a,-2L_a)\in \Tder(A)$ implies that $a\in\Nuc_m(A)$
and, therefore, $a\in\cen(A)$.
\end{proof}


\section{Nonexistence of ideals of finite codimension}
\label{finite.sec}

In this section $(V,[\,,\,,\,])$ will be a Lie triple system and
$U(V)$ --- the non-associative universal enveloping algebra of $V$.
For any $a,b,c \in V$ we have
\begin{equation*}[L_a,L_b](c) = a(bc) - b(ac)=[a,b,c] \quad \hbox{and} \quad [a,b] = 0
\end{equation*}
in $U(V)$. The map $[L_a,L_b]$ is a derivation of $U(V)$ and
$L_{ax+xa} = L_a L_x + L_x L_a $ for any $a\in V, x \in U(V)$ by
(\ref{leftalternative.eq}).

Let $A$ be a finite-dimensional unital algebra and $\LNalt(A)$
--- its left generalized alternative nucleus. We are interested in
the existence of monomorphisms of L.t.s.
\begin{equation}
\label{embedding.eq} \iota \colon V \rightarrow \LNalt(A)
\end{equation}
such that $\iota(a)\iota(b) =\iota(b)\iota(a)$ for any $a,b \in V$.
By the universal property of $U(V)$ such a map induces a
homomorphism $\varphi \colon U(V) \rightarrow A$. The kernel of
$\varphi$ is an ideal of finite codimension whose intersection with
$V$ is trivial.

\medskip

Let $S_2$ be the two-dimensional simple L.t.s.\ generated by $e,f$
with
\begin{equation} \label{twodimensional.eq} [e,f,e] = 2e \quad
\hbox{ and } \quad [e,f,f] = -2f.
\end{equation}
\begin{lem}
\label{commutatorS2.lem} With $e,f$ as above,
\begin{equation*}
[e^n,f] = n(n-1) e^{n-1}.
\end{equation*}
holds in $U(S_2)$.
\end{lem}
\begin{proof}
Observe that $fe^n = f(e e^{n-1}) = e(f e^{n-1}) -
[L_e,L_f](e^{n-1}) = e(fe^{n-1}) - 2(n-1) e^{n-1}$. Repeating with
$fe^{n-1}$ we obtain
\begin{eqnarray*}
fe^n &=& e^nf - 2((n-1) + (n-2) + \cdots + 1) e^{n-1}\\ &=& e^nf
-n(n-1) e^{n-1},
\end{eqnarray*}
\end{proof}

Any semisimple L.t.s. contains a copy of $S_2$, see \cite{Lis}.
(This may be compared to the fact that any semisimple Lie algebra
contains a copy of $sl_2$.)

\begin{prop}
\label{semisimplefinite.prop} If $(V,[\,,\,,\,])$ is a semisimple
L.t.s., then the only proper ideal of $U(V)$ that has finite
codimension is the augmentation ideal $\ker \epsilon$.
\end{prop}
\begin{proof}
Given a proper ideal $I$ of $U(V)$ whose codimension is finite, the
set $V_0 = I \cap V$ is an ideal of the L.t.s.\ $V$. Therefore,
there exists another ideal $V_1$ with $V=V_0 \oplus V_1$ (see
\cite{Lis}). Both $V_0$ and $V_1$ are semisimple L.t.s., so either
$V_1 = 0$, or there exists a subsystem $\spann \langle e,f\rangle
\subseteq V_1$ with multiplication as in (\ref{twodimensional.eq}).
In the first case we have that $\ker \epsilon$, the ideal generated
by $V$, is contained inside $I$ and, hence, since the codimension of
$\ker \epsilon$ is 1, they are equal.

Assume now that we are in the second case. Since any
finite--codimensional proper ideal $I$ of $U(V)$ contains an element
of the form $p(e) = \alpha_0 1 + \alpha_1 e + \cdots + \alpha_{n-1}
e^{n-1} + e^n$ with $n >1$, then, by Lemma \ref{commutatorS2.lem},
it also contains $[[p(e),f],f],\dots],f] = n! (n-1)! e$. Therefore,
$e \in I$ which, by definition of $V_1$, is not possible.
\end{proof}

Proposition \ref{semisimplefinite.prop} shows that embeddings of the
type (\ref{embedding.eq}) do not exist for semisimple L.t.s. Since
any L.t.s.\ decomposes (as a vector space) as the direct sum  of a
semisimple subsystem and a solvable ideal (see \cite{Lis}), it is
clear that such embedding might only exist for solvable L.t.s. We
shall prove that, in fact, $V$ must be nilpotent.

\medskip

Let us denote the map $c\mapsto [a,b,c]$ by $D_{a,b}$. The vector
space $\L_S(V) = \spann \langle D_{a,b}\,\vert\, a, b \in V\rangle
\oplus V$ is a Lie algebra (see \cite{J}) with the bracket
\begin{equation}
[a,b] = D_{a,b} \quad \hbox{and} \quad [D_{a,b},c] = [a,b,c].
\end{equation}
This Lie algebra is called sometimes {\em the Lie envelope} of $V$.
It is ${\Z}_2$--graded with even part $\L_S^+(V) = \spann \langle
D_{a,b}\vert\, a,b\in V\rangle$ and odd part $\L_S^-(V) = V$.

Given any unital algebra $A$ generated, as a unital algebra, by a
subsystem $V$ of $\LNalt(A)$ with $[a,b]=0$ for any $a,b\in V$, we
shall often consider the Lie algebra $\L(V)$ generated by $\{
L_a\,\vert\, a \in V\}$. Usually, no explicit mention of $A$ will be
needed. Since $[L_a,L_b]$ is a derivation of $A$ (see the proof of
Lemma \ref{derivation.lem}) and $A$ is generated by $V$, it follows
that $\L(V) = \spann \langle [L_a,L_b]\,\vert\, a,b \in V\rangle
\oplus \spann \langle L_a \,\vert\, a \in V\rangle$. The algebra
$\L(V)$ is isomorphic to $\L_S(V)$ by $a\mapsto L_a$ and
$D_{a,b}\mapsto [L_a,L_b]$.

\medskip

It is a simple exercise to check that over algebraically closed
fields of characteristic zero, the only solvable non-nilpotent
two-dimensional L.t.s.\ is $R_2 = Fa\oplus Fb$ with
\begin{equation}
[a,b,a] = -b \quad \hbox{and} \quad [a,b,b] = 0.
\end{equation}

\begin{lem}
\label{solvable.lem} Let $V$ be a solvable non-nilpotent L.t.s. Then
there exists a homomorphic image of $V$ which contains a subsystem
isomorphic to $R_2$.
\end{lem}
\begin{proof}
For $V$ a solvable L.t.s, $\L_S(V)$ is a solvable Lie algebra (see
\cite{Lis}). The solvability of $\L_S(V)$ implies that there exists
a non-zero $v \in \L_S(V)$ and a homomorphism of Lie algebras
$\lambda \colon \L_S(V) \rightarrow F$ such that
\begin{equation}\label{eigenvector.eq}
[x,v] =\lambda(x) v
\end{equation}
for any $x \in \L_S(V)$. Observe that $\L^+_S(V) \subseteq
[\L_S(V),\L_S(V)]$ and hence $\lambda(\L^+_S(V)) = 0$.

Write $v$ as a sum of its even and odd components: $v=D+b$ with
$D\in\L^+_S(V)$ and $b\in \L^-_S(V)$. The odd part of the
identity~(\ref{eigenvector.eq}) with $x\in \L_S^+(V)$ implies that
$[V,V,b] = 0$. Setting $x = a \in V$ in (\ref{eigenvector.eq}) gives
\[D_{a,b} = \lambda(a) D\] as the even part, and  \[D(a) =
-\lambda(a)b\] as the odd part.

Assume that $\lambda$ is not identically equal to zero. Then we can
choose $a \in V$ with $\lambda(a) = 1$. For such $a$ we have that
$D=D_{a,b}$ and $D(a) =-b$ so $[a,b,a] = -b$. Since $[V,V,b] = 0$,
the subspace $\spann\langle a,b\rangle$ is a subsystem of $V$
isomorphic to $R_2$.

Now, if $\lambda$ happens to be identically equal to zero, it
follows that $D=0$ and $b\neq 0$ (since $v$ is non-zero), and that
$[V,b,V] = 0$. Hence, the one-dimensional subspace $\spann\langle b
\rangle$ is contained in the centre of $V$. The L.t.s.\
$V/\spann\langle b \rangle$ is solvable non-nilpotent (see
\cite{Hop,Lis}) and its dimension is lower than the dimension of
$V$. The result in this case can be obtained by induction.
\end{proof}

\begin{prop}
Given a non-nilpotent L.t.s.\ $V$ and an ideal $I$ of finite
codimension in $U(V)$, the intersection $I\cap V$ is non-zero.
\end{prop}
\begin{proof}
Without loss of generality we may assume that $V$ is solvable. By
Lemma~\ref{solvable.lem}  there exists an ideal $V_0$ and elements
$a,b \in V$ such that $V_0 \oplus \spann\langle a,b \rangle$ is a
subsystem of $V$ with $[a,b,a] \equiv -b \mod{V_0}$ and $[a,b,b]
\equiv 0 \mod{V_0}$.

By (\ref{precession.eq}), in $U(V)$ we have $x(yz) = \sum
x_{(1)}y_{(1)} \cdot \delta_{x_{(2)}, y_{(2)}}(z)$. With $x = a^n$,
$y = c\in V$ and $z = a$ we obtain
\begin{equation*}
a^n(ca) = \left\{
\begin{array} {l}
a^{n+1} c = a\cdot a^n c \\
\\
a^nc \cdot a + \sum (a^n)_{(1)} \delta_{(a^n)_{(2)},c}(a) \\
\enskip\enskip \equiv a^n c \cdot a + n a^{n-1} \delta_{a,c}(a)
\mod{U(V)_{n-1}}
\end{array}
\right.
\end{equation*}
where the last congruence follows from (\ref{zerodegree.eq}). Hence
$[a^nc,a] \equiv -n a^{n-1} \delta_{a,c}(a) \equiv -\frac{n}{2}
a^{n-1}[a,c,a] \mod{U(V)_{n-1}}$. After $n$ commutations we get
\begin{equation*}
[\dots [[a^nc,a],a],\dots,a] =(-1)^n
\frac{n!}{2^n}[a,[a,[\dots,[a,c,a],\dots],a],a]
\end{equation*}
where we have replaced the congruence modulo $U(V)_0 = F$ by the
equality because both sides lie inside $\ker \epsilon$. In the
particular case $c = b$ we have $[\dots[[a^nb,a],a],\dots,a] =
\frac{n!}{2^n}(b+v_0)$ with $v_0 \in V_0$.

Any finite-codimensional ideal $I$ contains an element of the form
$p(a) = \alpha_0 1 + \alpha_1 a + \cdots + \alpha_{n-1} a^{n-1} +
a^n$. It also contains $p(a) b$ and $[\dots[p(a)b,a],\dots,a]$ where
the commutator is taken $n$ times. Therefore, $I$ also contains the
nonzero element $\frac{n!}{2^n}(b+ v_0)$.
\end{proof}

We have seen that faithful representations of the type
(\ref{embedding.eq}) can only exist for nilpotent L.t.s. It turns
out that for nilpotent L.t.s. these representations, if exist, have
very specific structure. Name, assuming that in (\ref{embedding.eq})
the algebra $A$ is generated by $\iota(V)$, we shall prove that
there exists a nilpotent ideal $R$ such that $A/R$ is a commutative
associative algebra over $F$ with no nontrivial nilpotent elements.
First, we need some lemmas.

\begin{lem}
Let $A$ be a finite-dimensional unital algebra, $a \in \LNalt(A)$
and $L_a = (L_a)_s + (L_a)_n$ --- the Jordan--Chevalley
decomposition of $L_a$ in $\Endo(A)$. Then there exist $a_s, a_n \in
\LNalt(A)$, the semisimple and nilpotent parts of $a$, with $(L_a)_s
= L_{a_s}$ and $(L_a)_n = L_{a_n}$.
\end{lem}
\begin{proof}
Recall that given $(d,d',d'') \in \Tder(A)$, its semisimple and
nilpotent parts can be calculated componentwise:
$(d,d',d'')_s=(d_s,d'_s,d''_s)$ and $(d,d',d'')_n=
(d_n,d'_n,d''_n)$, where both $(d_s,d'_s,d''_s)$ and
$(d_n,d'_n,d''_n)$ are also ternary derivations. Recall also that
$(d,d',-d) \in \Tder(A)$ if and only if $d=L_a$ and $d'=T_a$ with $a
\in \LNalt(A)$. Now, for any $a \in \LNalt(A)$ we have that
$((L_a)_s,(T_a)_s,-(L_a)_s)$ and $((L_a)_n,(T_a)_n,-(L_a)_n) \in
\Tder(A)$, which implies that $(L_a)_s = L_{a_s}$ and $(L_a)_n =
L_{a_n}$ for some $a_s, a_n \in \LNalt(A)$.
\end{proof}

Let us complete $V$ inside $A$ by adding the semisimple and
nilpotent parts of all its elements; it turns out that such
completion retains some fundamental properties of $V$:

\begin{lem}
\label{JordanChevalley.lem}
 Let $A$ be a finite dimensional unital algebra. Given any subsystem $V
\leq \LNalt(A)$ such that
\begin{itemize}
\item[i)] $V$ generates $A$ as a unital algebra,
\item[ii)] $[a,b] = 0$ for all $a,b \in V$,
\item[iii)] $V$ is nilpotent,
\end{itemize}
there exists in $\LNalt(A)$ a subsystem $\hat{V}$ containing $V$ and
satisfying i), ii) and iii), and such that $a_s, a_n \in \hat{V}$
for any $a \in \hat{V}$. Moreover, $a_s \in \cen(A)$ for any $a\in
\hat{V}$ and $\{ a_n\,\vert\, a\in \hat{V}\}$ is an ideal of
$\hat{V}$.
\end{lem}
\begin{proof}
Since $V$ generates $A$ and $[a,b] = 0$ for any $a,b \in V$, the Lie
algebra $\L(V)$, generated by $\{ L_a \,\vert\, a \in V\}$ is
isomorphic to $\L_S(V)$. By \cite{Hop} the latter algebra, and hence
the former, is nilpotent.

 By the properties of the
Jordan--Chevalley decomposition (see \cite{Hump}) $(\ad_{L_a})_s =
\ad_{L_{a_s}}$ and $(\ad_{L_a})_n = \ad_{L_{a_n}}$. The operators
$\ad_{L_{a_s}}$ and $\ad_{L_{a_n}}$ can be expressed as polynomials
in $\ad_{L_a}$ with zero constant term. In particular,
$\ad_{L_{a_s}}$ leaves $\L(V)$ stable with a nilpotent action. By
the semisimplicity of $\ad_{L_{a_s}}$ this means that
$[L_{a_s},\L(V)] = 0$. Hence $a_s \in \cen(A)$ by
Lemma~\ref{center.lem}.

As $\L(V)$ is nilpotent, there exists a basis of $A$ where $\L(V)$
is represented by upper triangular matrices. Hence, for any $a,b \in
V$ the operator $L_{a_s + b_s}$ is semisimple, while $L_{a_n + b_n}$
is nilpotent. Moreover, $a_s + b_s \in \cen(A)$ implies that
$[L_{a_s + b_s}, L_{a_n + b_n}] = 0$. By the uniqueness of the
Jordan--Chevalley decomposition we obtain that $(L_{a+b})_s =
L_{a_s} + L_{b_s}$ and $(L_{a+b})_n = L_{a_n} + L_{b_n}$. In
particular, $(a+b)_s = a_s + b_s$ and $(a+b)_n = a_n + b_n$.

Let $\hat{V} = \{ a_s + b_n \,\vert\, a,b \in V\}$. By the previous,
$\hat{V}$ is a vector subspace of $\LNalt(A)$ and, since $(a_s +
b_n)_s = a_s$ and $(a_s + b_n)_n = b_n$, $\hat{V}$ contains the
semisimple and nilpotent components of its elements. We also know
that $a_s \in\cen(A)$ for any $a\in \hat{V}$.

Given $a,a',a'',b,b',b'' \in V$ we have that
\begin{equation*}
[a_s+b_n,a'_s +b'_n] = [b_n,b'_n] = [b_s+b_n,b'_s+b'_n] = [b,b'] =
0,
\end{equation*}
so $\hat{V}$ satisfies ii). Moreover,
\begin{eqnarray*}
[a_s + b_n, a'_s + b'_n, a''_s+b''_n] &=& [L_{a_s+b_n},L_{a'_s +
b'_n}](a''_s + b''_n)\\ &=& [L_{b_n},L_{b'_n}](a''_s + b''_n) \\
&=&[b_n,b'_n,b''_n]\\ &=&[b,b',b'']
\end{eqnarray*}
implies that $\hat{V}$ is a subsystem of $\LNalt(A)$ and that
$[\hat{V},\hat{V},\hat{V}]\subseteq [V,V,V]$. In terms of the lower
central series for $\hat{V}$ and $V$ (see \cite{Hop}) this says that
$\hat{V}^1\subseteq V^1$. Assuming that $\hat{V}^n\subseteq V^n$, we
have $\hat{V}^{n+1}=[\hat{V}^n,\hat{V},\hat{V}] +
[\hat{V},\hat{V},\hat{V}^n]\subseteq [V^n,\hat{V},\hat{V}] +
[\hat{V},\hat{V},V^n]\subseteq [V^n,V,V]+[V,V,V^n]\subseteq
V^{n+1}$. The nilpotency of $\hat{V}$ follows from this observation
and the nilpotency of $V$.

Finally, the left multiplication operator by $[a,b,c]$ is obtained
as the commutator $[[L_a,L_b],L_c]$; in an adequate basis of $A$ it
is represented as a commutator of upper triangular matrices.
Therefore, it is nilpotent and $[a,b,c] = [a,b,c]_n$. Since
$[\hat{V},\hat{V},\hat{V}]\subseteq [V,V,V]$ it follows that
$[\hat{V},\hat{V},\hat{V}]\subseteq \{ a_n\,\vert\, a\in \hat{V}\}$.
In particular, the latter set is an ideal of $\hat{V}$.
\end{proof}

\begin{lem}
\label{basicnilpotency.lem} Let $A$ be a finite-dimensional unital
algebra and let $V$ be a subsystem of $\LNalt(A)$. Assume that
\begin{itemize}
\item[i)] $a = a_n$ for any $a\in V$,
\item[ii)] $[a,b]=0$ for all $a,b\in V$.
\end{itemize}
Then the subalgebra generated by $V$ is nilpotent.
\end{lem}
\begin{proof}
Assume, as before, that $A$ is generated by $V$ as a unital algebra.

There exists an element of $V$ that lies in the centre of $A$.
Indeed, the nilpotency of $V$ implies that $\L(V)$ consists of
nilpotent transformations \cite{Hop}, which, in turn, implies that
the centre of $\L(V)$ is non-zero. Given $0\neq D + L_a \in
\cen(\L(V))$ with $D\in \L^+(V)$, for any $b\in V$ the equality $0 =
[D+ L_a , L_b] = L_{D(b)} +[L_a,L_b]$ implies that $D = 0$ and
$[L_a,L_b] = 0$. Therefore $0\neq a \in \cen(A)$ by Lemma
\ref{center.lem}.

We shall use induction on the dimension of $V$. The case $\dim V =
0$ is obvious. Given $V$ with $\dim V = n+1$, choose $0\neq a\in
Z(A)\cap V$ as above and consider the ideal $aA$. The quotient
algebra $A/aA$ is generated, as a unital algebra, by the quotient
$(V+ aA)/aA$ of $V$. Thus we can apply the hypothesis of induction
to conclude that $\alg\langle V+aA/aA\rangle = \alg\langle V\rangle
/aA$ is nilpotent.

Let us denote the ideal $\alg\langle V\rangle$ by $A_0$, and the
linear span of all products of $N$ elements of $A_0$, regardless of
the order of the parentheses, by $A_0^N$. From the nilpotency of
$\alg\langle V\rangle /aA$ we deduce that there exists $N$ such that
$A_0^N \subseteq a A$. Moreover, any product involving $2N$ elements
of $A_0$ lies in the ideal $aA_0$, since is of the form $u_1u_2$
where at least one of the factors involves at least $N$ elements
and, therefore, lies in $A_0^N\subseteq aA$, and the other factor
belongs to $A_0$.

Let us fix $N$ such that $A_0^N\subseteq aA_0$ and prove by
induction that $A_0^{N^k}\subseteq a^k A_0$. Since $a\in V $ is
nilpotent, this will imply that $A_0$ is nilpotent, as desired.
Assume that $A_0^{N^{k-1}}\subseteq a^{k-1} A_0$. Any product of
$N^k$ elements in $A_0$ can be written as a product of $N$ factors,
each belonging to $A_0$, and at least one of them lying in
$A_0^{N^{k-1}}\subseteq a^{k-1} A_0$. Since $a^{k-1}$ is in the
centre of $A$, the whole product lies in $a^{k-1}A_0^N\subseteq
a^kA_0$.
\end{proof}

Finally, we are in the position to prove Theorem~\ref{main.thm}.

\begin{proof}[Proof of Theorem~\ref{main.thm}]
By Lemma \ref{JordanChevalley.lem} we can assume that $V$ contains
the semisimple and nilpotent components of all its elements. Let $Q
=\alg\langle a_s\,\vert\, a\in V\rangle \subseteq \cen(A)$ and let
$R$ be the ideal generated by $\{ a_n\,\vert\, a\in V\}$. Clearly $A
= Q + R$.

For any nilpotent element $x\in Q$, $L_x$ belongs to $\alg_1\langle
L_{a_s}\,\vert\, a\in V\rangle$. This algebra is abelian and all its
elements are semisimple transformations. But $x\in\cen(A)$ implies
that $L_x$ is nilpotent so $L_x = 0$ and $x$ must be zero. Hence $Q$
is a commutative associative finite dimensional algebra without
nonzero nilpotent elements.

Since $a_s\in \cen(A)$, it follows that $A = Q\alg_1\langle
a_n\,\vert\, a\in V\rangle$. We can apply Lemma
\ref{basicnilpotency.lem} to the algebra $\alg_1\langle a_n\,\vert\,
a\in V\rangle$ and the subsystem $\{ a_n\,\vert\, a\in V\}$ to
conclude that $\alg\langle a_n\,\vert\, a\in V\rangle$ is nilpotent.
The ideal $R$ decomposes as $R =Q\alg\langle a_n\,\vert\, a\in
V\rangle$, so it is also nilpotent. Its nilpotency implies that
$Q\cap R = 0$.
\end{proof}


\section{The universal enveloping algebras of a L.t.s.\ and its Lie envelope}
\label{construction.sec}

The following construction is based on the known construction of a
Bruck loop starting from a group whose every element has a square
root. Namely, any such group with the product $g*h=g^{\frac{1}{2}} h
g^{\frac{1}{2} }$ becomes a Bruck loop. Observe that the
linearisation of the identity $g = r(g) r(g)$ with $r(g) =
g^{\frac{1}{2}}$ in an H--bialgebra reads as $x = \sum r(x_{(1)})
r(x_{(2)})$ for some map $r$.

\medskip

Let $L$ be a Lie algebra over a field $F$ of characteristic $\neq
2$.

\begin{lem}
The linear map $q\colon U(L) \rightarrow U(L)$ defined by $x\mapsto
\sum x_{(1)}x_{(2)}$ is bijective.
\end{lem}
\begin{proof}
Consider the Poincar\'e--Birkhoff--Witt filtration $U(L) =
\bigcup_{n\geq 0} U_n$ of $U(L)$. Given $a_1,\dots, a_n \in L$,
\[q(a_1\cdots a_n) \equiv 2^n a_1 \cdots a_n \mod{U_{n-1}}.\]
Since $q$ preserves the filtration, it follows that it is bijective
on each $U_n$.
\end{proof}

Let $r$ be the inverse of $q$. Clearly, for any $x \in U(L)$ we have
that $x = \sum r(x)_{(1)} r(x)_{(2)}$. Furthermore, $q$ being a
coalgebra isomorphims implies that $r$ is also a coalgebra
isomorphism. Therefore,
\[
x = \sum r(x_{(1)})r(x_{(2)})
\]

The product on $U(L)$ can be modified with the help of the map $r$
as follows:
\[
x*y = \sum r(x_{(1)}) y r(x_{(2)}).
\]
With this product $U(L)$ becomes a unital non-associative algebra.
In fact, since $r$ is a homomorphism of coalgebras,  $U(L)$ carries
the structure of an H--bialgebra.

\begin{lem}\label{identity.lem} For all $x,y$ in $U(L)$
\[\sum r(x_{(1)}*(y*x_{(2)})) = \sum r(x_{(1)}) r(y) r(x_{(2)}).\]
\end{lem}
\begin{proof} Indeed,
\begin{eqnarray*}
\sum x_{(1)} *(y*x_{(2)}) &=&\sum r(x_{(1)})r(y_{(1)})x_{(3)}r(y_{(2)})r(x_{(2)}) \\
&=& \sum r(x_{(1)})r(y_{(1)})r(x_{(2)})r(x_{(3)})r(y_{(2)})r(x_{(4)}) \\
&=& \sum \left( r(x_{(1)})r(y)r(x_{(2)})\right)_{(1)} \left(
r(x_{(1)})r(y)r(x_{(2)})\right)_{(2)}
\end{eqnarray*}
which proves the lemma.
\end{proof}

\begin{prop}\label{mult.prop}
The algebra $(U(L),*)$ satisfies
\begin{itemize}
\item[i)] $\sum x_{(1)} * (y *(x_{(2)}*z)) = \sum (x_{(1)}*(y*x_{(2)}))*z$.
\item[ii)] $a*b = b*a$ for any $a,b \in L$.
\item[iii)] $a*(b*c) - b*(a*c) = \frac{1}{4}[[a,b],c]$ for any $a,b,c \in L$.
\end{itemize}
\end{prop}
\begin{proof}
We shall only check part i); it follows from
Lemma~\ref{identity.lem} by
\begin{multline*}
\sum x_{(1)}*(y*(x_{(2)}*z)) = \sum r(x_{(1)}) r(y_{(1)})
r(x_{(2)})z r(x_{(3)}) r(y_{(2)})r(x_{(4)}) \\
=\sum r(x_{(1)}*(y*x_{(2)}))_{(1)}zr(x_{(1)}*(y*x_{(2)}))_{(2)} =
\sum (x_{(1)}*(y*x_{(2)}))*z.
\end{multline*}
\end{proof}

Given a L.t.s.\ with the product $[\,,\,,\,]$ and a scalar $\mu$,
the new product $[\,,\,,\,]' = \mu^2 [\,,\,,\,]$ also defines a
L.t.s.\ that is isomorphic to the original L.t.s.\ under $x\mapsto
\mu x$.

\begin{cor}
Let $V$ be a L.t.s.\ and $\L_S(V)$ --- the Lie envelope of $V$. The
unital subalgebra of $(U(\L_S(V)),*)$ generated by $V$ is isomorphic
to the universal enveloping algebra of $V$ considered as a Bol
algebra with the trivial binary product.
\end{cor}
\begin{proof}
Define $[a,b,c]' =\frac{1}{4}[a,b,c]$ and let $Q$ the subalgebra of
$(U(L),*)$ generated by $V$. The universal property of
$U(V,[\,,\,,\,]')$ together with Proposition~\ref{mult.prop} implies
that there exists an epimorphism from $U(V,[\,,\,,\,]')$ to $Q$.
Since $a_1*(\cdots *(a_{n-1}*a_n)) \equiv a_1\cdots a_n
\mod{U_{n-1}}$ with $a_1,\dots, a_n \in V$, it follows that $Q$
admits a PBW--type basis. The epimorphism from $U(V,[\,,\,,\,]')$ to
$Q$ maps the PBW basis of $U(V,[\,,\,,\,]')$ to this basis, so it is
an isomorphism. However, as $(V,[\,,\,,\,]')$ and $(V,[\,,\,,\,])$
are isomorphic, their universal enveloping algebras also are.
\end{proof}


\section{Ideals in the enveloping algebras of simple L.t.s.} \label{simple.sec}

\begin{lem}\label{proper.ideal.lem}
Assume that $V$ is a simple L.t.s.\ satisfying the following
condition: all elements of the universal enveloping algebra $U(V)$
that commute with $V$, are of the form $c+x$ where $c$ is a scalar
and $x$ is in $V$. Then the only proper ideal of $U(V)$ is the
augmentation ideal.
\end{lem}
\begin{proof}
Suppose the conditions of the lemma are satisfied. Let $I\subset
U(V)$ be an ideal, and take some $r\in I$. There exists an element
$x\in V$ such that $r'=rx-xr\neq 0$. It is clear that $r'\in I$ and
$\deg r'<\deg r$, where the degree is taken with respect to the PBW
filtration. Hence, $I$ necessarily contains a nonzero element $u$ of
degree at most 1. If $u$ is a scalar, then $I=U(V)$. If $\deg u =1$,
the space of all linear combinations of (possibly iterated) brackets
containing $u$, is an ideal of $V$ and, hence, coincides with $V$.
All these brackets are in $I$, therefore, $I$ contains $V$.
\end{proof}

If $a$ is an element of $V$ and $r\in U(V)_n$, the commutator
$ar-ra$ belongs to $U(V)_{n-1}$. In fact, it is possible to write an
explicit formula for the terms of degree $n-1$ in this commutator.
\begin{lem}
Let $\{x_{k}\}$ be a basis for $V$ and $r\in U(V)$ --- a monomial in
the $x_k$. Then
\begin{equation}\label{partial.derivatives.eq}
ar-ra=\frac{1}{2}\sum_{i,j}\ [a,x_i,x_j]\cdot \ddx{i}\ddx{j}r+
\text{lower degree terms}.
\end{equation}
\end{lem}
Here the partial derivative $\partial/\partial x_i$ of a
non-associative monomial is defined by setting $\partial/\partial
x_i(uv)=u\partial/\partial x_i(v)+\partial/\partial x_i(u)v$ with
$\partial/\partial x_i(x_j)=1$ if $i=j$ and 0 otherwise.

\begin{proof}
The vector space $U(V)_n/U(V)_{n-1}$ is spanned by classes of
elements of the form $b^n$ with $b\in V$, so it is sufficient to
verify (\ref{partial.derivatives.eq}) for $p=b^n$.

Modulo terms of degree $n-2$ and smaller we have
\begin{align*}
ab^n-b^na&= a(bb^{n-1})-b^na \\
&=[L_a,L_b](b^{n-1})+b(ab^{n-1})-b^na\\
&=\sum_{i+j=n-2} b^i ([L_a,L_b](b)\cdot b^j) +b(ab^{n-1})-b^na\\
&=(n-1)[L_a,L_b](b)\cdot b^{n-2}+b(ab^{n-1}-b^{n-1}a)\\
&=(n-1)[L_a,L_b](b)\cdot b^{n-2}+(n-2)[L_a,L_b](b)\cdot b^{n-2}+\ldots+[L_a,L_b](b)\cdot b^{n-2}\\
&=\frac{n(n-1)}{2}[L_a,L_b](b)\cdot b^{n-2}.
\end{align*}
The last expression coincides with the right-hand side of
(\ref{partial.derivatives.eq}).
\end{proof}

\medskip

Let $x,y,z$ be a set of generators for the Lie algebra $so(3)$ with
$[x,y]=z$, $[y,z]=x$ and $[z,x]=y$. We shall consider $so(3)$ as a
simple L.t.s.\ by setting $[a,b,c]=[[a,b],c]$. Let $\widetilde{S}_2$
be the 2-dimensional subsystem spanned by $x$ and $y$. Over the
complex numbers $\widetilde{S}_2$ is isomorphic to the L.t.s.\ $S_2$
mentioned in Section~\ref{stuff.sec}; the isomorphism is given by
$e=-x+y\sqrt{-1}$, $f=x+y\sqrt{-1}$.
\begin{prop}
Both $so(3)$ and $\widetilde{S}_2$ satisfy
Conjecture~\ref{main.conj}.
\end{prop}
\begin{proof}
The products of the form $z^n(x^py^q)$ with $n,p,q$ non-negative
integers, form a basis for the universal enveloping algebra of
$so(3)$ considered as a Lie triple system. In our case,
(\ref{partial.derivatives.eq}) reads as
\begin{multline*}
z^n(x^py^q)\cdot z- z\cdot z^n(x^py^q)\\
=-\frac{n(p+q)}{2}z^{n-1}(x^py^q)+\frac{p(p-1)}{2}z^{n+1}(x^{p-2}y^q)+
\frac{q(q-1)}{2}z^{n+1}(x^py^{q-2})+\ldots
\end{multline*}
where the omitted terms are of degree $n+p+q-2$ and smaller.
Similarly,
\begin{multline*}
z^n(x^py^q)\cdot x- x\cdot z^n(x^py^q)\\
=\frac{n(n-1)}{2}z^{n-2}(x^{p+1}y^q)-\frac{p(n+q)}{2}z^{n}(x^{p-1}y^q)+
\frac{q(q-1)}{2}z^{n}(x^{p+1}y^{q-2})+\ldots
\end{multline*}
and
\begin{multline*}
z^n(x^py^q)\cdot y- y\cdot z^n(x^py^q)\\
=\frac{n(n-1)}{2}z^{n-2}(x^py^{q+1})+\frac{p(p-1)}{2}z^{n}(x^{p-2}y^{q+1})-
\frac{q(n+p)}{2}z^{n}(x^py^{q-1})+\ldots
\end{multline*}
Now, suppose that there exists an element $r$ of the universal
enveloping algebra of $so(3)$ considered as a Lie triple system, of
degree $N>1$, which commutes with $x,y$ and $z$. This element has
the form
\[ r=\sum_{n+p+q=N} \alpha_{n,p,q} z^n(x^py^q) +\text{lower\ degree\
terms.}\] The requirement that $rz-zr$ has no terms of degree $N-1$
imposes linear conditions on the coefficients $\alpha_{n,p,q}$,
similar conditions come from $rx-xr$ and $ry-yr$. Explicitly, these
conditions are as follows:
\begin{align*}
-(p+q)(n+2)\alpha_{n+2,p,q}&+(p+1)(p+2)\alpha_{n,p+2,q}+(q+1)(q+2)\alpha_{n,p,q+2}=0,\\
(n+1)(n+2)\alpha_{n+2,p,q}&-(n+q)(p+2)\alpha_{n,p+2,q}+(q+1)(q+2)\alpha_{n,p,q+2}=0,\\
(n+1)(n+2)\alpha_{n+2,p,q}&+(p+1)(p+2)\alpha_{n,p+2,q}-(n+p)(q+2)\alpha_{n,p,q+2}=0.
\end{align*}
The determinant of the corresponding $3\times 3$-matrix is equal to
$2(n+2)(p+2)(q+2)(n+p+q+1)^2$ and it follows that all the
$\alpha_{n,p,q}$ are zero and, hence, $\deg{r}<N$, which gives a
contradiction.

The argument for $\widetilde{S}_2$ is entirely similar.
\end{proof}

Let $(\, ,)$ be a non-degenerate symmetric bilinear form on a vector
space $V$ of dimension greater than $1$. Define a ternary bracket on
$V$ by
\[[a,b,c]=(a,c)b-(b,c)a.\] A straightforward verification
shows that $V$ with this bracket satisfies all the axioms of a Lie
triple system.

If $I$ is an ideal in $V$, $[I,V,V]\subseteq I$, that is,
$(v,x)u-(v,u)x\in I$ for any $x\in I$ and $u,v\in V$. Hence,
$(v,x)u\in I$ for any $x\in I$ and this means that $I$ is either
trivial, or coincides with $V$. Therefore, $V$ is simple.

\begin{prop}
The L.t.s.\ $V$ satisfies Conjecture~\ref{main.conj}.
\end{prop}
\begin{proof}
Fix a basis $\{x_{k}\}$ for $V$, $n\geq k\geq 1$, and let $r\in
U(V)$ be homogeneous of degree greater than 1. The condition  $x_k
r-r x_k=0$ implies, by (\ref{partial.derivatives.eq}), that
\[\sum_{i,j} [x_k, x_i, x_j]\ddx{i}\ddx{j}r=0,\]
that is,
\[\sum_{i,j}((x_k,x_j)x_i-(x_i,x_j) x_k)\ddx{i}\ddx{j}r=0.\]
Assuming that the basis $\{x_{k}\}$ is orthonormal, we get
\[x_k  \sum_{i} \frac{\partial^2}{\partial x_i^2} r=
\sum_{i} x_i\ddx{i}\ddx{k}r = (m-1)\ddx{k}r,\] where $m=\deg r$. If
$\ddx{k}r=0$ for some $k$ it follows that $\sum_{i}
\frac{\partial^2}{\partial x_i^2}r=0$ and, hence, that $\ddx{k}r=0$
for all $k$. In this case $r$ is a constant, so we can assume that
$\ddx{k}r\neq 0$ for all $k$ and that $\sum_{i}
\frac{\partial^2}{\partial x_i^2}r\neq 0$.

Let us write $\psi$ for $\sum_{i} \frac{\partial^2}{\partial
x_i^2}r$. We have
\begin{equation}\label{psi.eq}(m-1)\ddx{k} r= x_k\psi\end{equation}
and, hence,
\[(m-1)x_k\ddx{k} r= x_k^2\psi\]
and
\[m(m-1) r= q\psi\]
with $q=\sum_{i=1}^n x_i^2$. It follows that
\[m(m-1)\ddx{k} r= 2 x_k\psi+ q\ddx{k}\psi\]
which implies
\[(m-2)x_k\psi = q\ddx{k}\psi.\]
If $m=2$ this means that $r$ is a scalar multiple of $q$. If $m\neq
2$ we have that $\psi=q\psi_0$ with $\psi_0\neq 0$, and $m(m-1) r=
q^2\psi_0$. It is readily seen that $\psi_0$ satisfies
\[(m-4)x_k\psi_0 = q\ddx{k}\psi_0.\]
If $m=4$ this implies that $r$ is a scalar multiple of $q^2$;
otherwise the above manipulations can be repeated. Eventually, this
process has to stop and in the end we get that $m=2l$ and that, up
to a multiplication by a scalar, $r=q^l$.

Now, (\ref{psi.eq}) can be re-written as
\[(2l-1)\cdot 2x_k l q^{l-1}= x_k (2n l q^{l-1}+4l(l-1)q^{l-1}).\]
This gives $n=1$ and it follows that $x_k r-r x_k=0$ cannot be
satisfied for all $k$.
\end{proof}



{\small \begin{thebibliography}{99}


\bibitem{Hop}
N. C. Hopkins, Nilpotent Ideas in Lie and Anti--Lie Triple Systems,
J. Algebra {\bf 178} (1995), 480--492.

\bibitem{Hump}
J. E. Humphreys, \emph{Introduction to Lie Algebras and
Representation Theory}\/, Springer Verlag, New York, 1972.


\bibitem{J} N. Jacobson, General representation theory of Jordan algebras, Trans. Amer. Math. Soc. {\bf 70} (1951), 509--530.

\bibitem{Ki}
H. Kiechle, \emph{Theory of $K$--loops}\/ Lecture Notes in
Mathematics, 1778. Springer--Verlag, Berlin, 2002.

\bibitem{Lis}
W. G. Lister, A structure theory of Lie triple systems, \emph{Trans.
Amer. Math. Soc.} {\bf 72} (1952), 217--242.

\bibitem{P} J.M. P\'erez--Izquierdo,
An Envelope for Bol Algebras, {\em J. Algebra\/} {\bf 284} (2005),
480--493.

\bibitem{P2} J.M. P\'erez--Izquierdo,
Algebras, hyperalgebras,  nonassociative bialgebras  and loops,
http://mathematik.uibk.ac.at/\hskip 0.02pt mathematik/\hskip 0.01pt
jordan/.

\bibitem{PS} J.M. P\'erez--Izquierdo and I.P. Shestakov,
An Envelope for Malcev Algebras, {\em J. Algebra\/} {\bf 272}
(2004), 379--393.

\bibitem{S} L.V. Sabinin, {\em Smooth Quasigroups and Loops\/} (Mathematics
and Its Applications, 492, Kluwer Academic Publishers, 1999).

\bibitem{SM-Orange} P.\ Miheev\ and\ L.\ Sabinin,
{\em Quasigroups and differential geometry}, Quasigroups and loops:
theory and applications, 357--430, Heldermann, Berlin, 1990.

\bibitem{SU} I.P. Shestakov and U.U. Umirbaev, Free Akivis algebras,
primitive elements, and hyperalgebras,  {\em J. Algebra\/} {\bf 250}
(2002), no. 2, 533--548.

\bibitem{Sw} E.M. Sweedler, \emph{Hopf algebras}\/ Mathematics Lecture
Note Series, W. A. Benjamin, Inc., New York, 1969.


\end{thebibliography}}
 \end{document}